\documentclass{article}
\usepackage{amssymb,graphpap,xcolor,xypic,color}

\usepackage{geometry}
\usepackage{doc}

\usepackage{url}

\usepackage{graphicx}
\usepackage{epstopdf}
\title{Geometry of Regular Algebras of Global Dimension 4 related to Graded Skew Clifford Algebras of Global Dimension 3}
\author{Manizheh Nafari}
\date{}

\begin{document}

% Add the title section.
\maketitle

% Add an abstract.
\abstract{
\noindent We compute point schemes of some regular algebras in \cite{SV2} using (Wolfram) Mathematica. These algebras are Ore extensions of regular graded skew Clifford algebras of global dimension 3 (c.f., [chapter 4, \cite{N1}]).}

% Add various lists on new pages.
\pagebreak
\tableofcontents
\listoffigures

%\pagebreak
%\listoftables

% Start the paper on a new page.
\pagebreak

%
% Body text.
%
\section{Introduction}
\label{introduction}

\noindent M. Artin, W. Schelter, J. Tate, and M. Van den Bergh introduced the notion of non-commutative regular algebras and invented new methods in algebraic geometry in the late 1980s to study them (\cite{AS}, \cite{ATV1}, \cite{ATV2}). Such algebras are viewed as non-commutative analogues of polynomial rings; indeed, polynomial rings are examples of regular algebras.\\

\noindent By the 1980s, a lot of algebras had arisen in quantum physics, specifically quantum groups, and many traditional algebraic techniques failed on these new algebras. In physics, quantum groups are viewed as algebras of non-commuting functions acting on some ``non-commutative space''(\cite{D}). In the early 1980s, E. K. Sklyanin, a physicist, constructed a family of graded algebras on four generators (\cite{S}). These algebras were later proved to depend on an elliptic curve and an automorphism (\cite{FO}). By the late 1980s, it was known that many of the algebras in quantum physics are regular algebras; in particular, the family of algebras constructed by Sklyanin consists of regular algebras.\\

\noindent In 2010, T. Cassidy and M. Vancliff introduced a class of algebras that provide an ``easy'' way to write down some quadratic regular algebras of global dimension $n$  where $n \in \mathbb{N}$ (\cite{CV}). In fact, they generalized the notion of a graded Clifford algebra and called it a graded skew Clifford algebra (see [Definition 3.1, \cite{N2}]).\\

\noindent In this lecture notes, we compute point schemes of some regular algebras in \cite{SV2} using (Wolfram) Mathematica. These algebras are Ore extensions of regular graded skew Clifford algebras of global dimension 3 (c.f., [chapter 4, \cite{N1}]).\\

\noindent Remark: We compute the coordinate $(\alpha_1,\alpha_2,\alpha_3,\alpha_4) \in \mathbb{P}^3$. The computation of $(\beta_1,\beta_2,\beta_3,\beta_4) \in \mathbb{P}^3$ will be added in the near future.

\section{Acknowledgement}

\noindent I would like to thank Michaela Vancliff (my Ph.D. advisor during my studies 2007-2011 at the University of Texas at Arlington) for her invaluable advice and for explaining essential concepts on non-commutative algebraic geometry.

\pagebreak

%\section{Background, Preliminary, and Related Work}

%\section{Main Content Sections}

%\subsection{Multiple Outline Levels}

%\subsection{Tables and Figures}

%\begin{table}
%\centering
%\begin{tabular}{|c|c|c|}\hline
%Column 1 & Column 2 & Column 3 \\\hline\hline
%a & b & c \\
%d & e & f \\
%g & h & i \\\hline
%\end{tabular}

%\caption{A sample table}
%\label{table-sample}
%\end{table}

%\begin{figure}
%\centering
%\includegraphics[width=2in]{space.jpg}

%\caption{A sample figure}
%\label{figure-sample}
%\end{figure}

%\section{Geometry of Graded Skew Clifford Algebras of Global Dimension 3}

%\subsection{Triangle}

%\subsection{Conic Union Line}

%We may as well include a second figure also, shown in Figure~\ref{figure-sample2}.  The same image file is used, but note how it can be resized.  Again, observe how the positions of the tables and figures do not necessarily match their positions in the source file, reiterating the aforementioned \LaTeX\ functionality for deciding where these items go in the final document.  You provide an approximate location, and \LaTeX\ does the rest.

%\begin{figure}
%\centering
%\includegraphics[width=1in]{space.jpg}

%\caption{Another sample figure}
%\label{figure-sample2}
%\end{figure}

%\subsection{Nodal Cubic Curve}

\section{Geometry of Graded Skew Clifford Algebras of Global Dimension 4}

\noindent Throughout this lecture notes, $\mathbb{K}$ denotes an algebraically closed field, $\text{char}(\mathbb{K}) \neq 2$, and $\mathbb{K}^{\times}$ denotes $\mathbb{K} \setminus \{ 0 \}$.

\subsection{Proposition 1:}

\noindent Suppose $q \in \mathbb{K}^{\times}$, and $q^ 2 \neq 1$. If the algebra
\[
A = \frac{\mathbb{K} \langle x_1,x_2,x_3,x_4 \rangle}{\langle g_1, \ldots, g_6 \rangle}
\]
where
\[
\begin{array}{ll}
g_1 = x_2 x_1 - q x_1 x_2, & \quad g_2 = x_2 x_3 - x_3 x_2,\\
g_3 =  x_3 x_1 - q x_1 x_3, & \quad g_4 = x_4 x_1 - x_1 x_4 - (q - q^{-1}) x_2 x_3,\\
g_5 = x_4 x_2 - q x_2 x_4, & \quad g_6 = x_4 x_3 - q x_3 x_4,
\end{array}
\]
, then $A$ has point scheme given by $\mathcal V(x_2,x_3) \cup \mathcal V(x_2 x_3 - x_1 x_4)$ (see Figure 1).

\begin{figure}[ht]
\centering
\includegraphics[width=2cm,height=2cm]{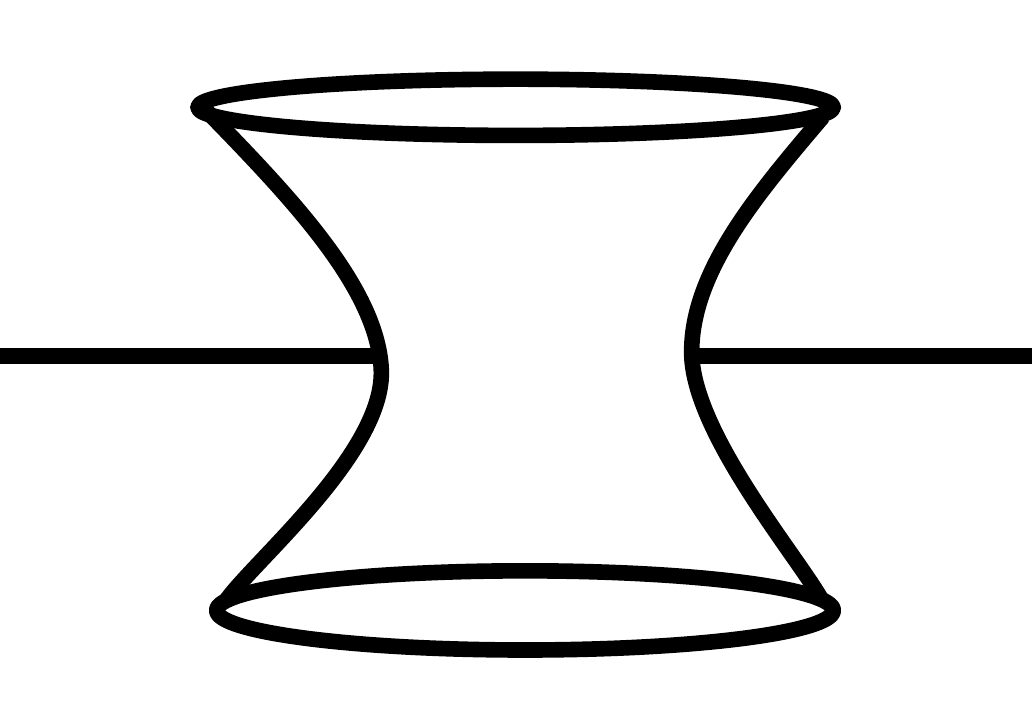}\\
\caption{Depiction of the Point Scheme in Proposition 1}
\end{figure}

\noindent Proof:\\

\noindent Suppose
\[
p = ((\alpha_1,\alpha_2,\alpha_3,\alpha_4),(\beta_1,\beta_2,\beta_3,\beta_4)) \in \mathbb{P}^3 \times \mathbb{P}^3.
\]
To find the point scheme $\mathcal{P}$ of $A$, we solve
\[
0 = g_1(p) = \alpha_2 \beta_1 - q \alpha_1 \beta_2,
\]
\[
  0 = g_2(p) = \alpha_2 \beta_3 - \alpha_3 \beta_2,
\]
\[
0 = g_3(p) = \alpha_3 \beta_1 - q \alpha_1 \beta_3,
\]
\[
0 = g_4(p) = \alpha_4 \beta_1 - \alpha_1 \beta_4 - (q - q^{-1}) \alpha_2 \beta_3,
\]
\[
0 = g_5(p) = \alpha_4 \beta_2 - q \alpha_2 \beta_4,
\]
\[
0 = g_6(p) = \alpha_4 \beta_3 - q \alpha_3 \beta_4,
\]
which yields $D E = F$, where
\[D =
\left[\begin{array}{cccc} \alpha_2 & -q \alpha_1 & 0 & 0 \\ 0 & - \alpha_3 & \alpha_2 & 0 \\ \alpha_3 & 0 & - q \alpha_1 & 0 \\ \alpha_4 & 0 & - (q - q^{-1}) \alpha_2 & - \alpha_1 \\ 0 & \alpha_4 & 0 & - q \alpha_2 \\ 0 & 0 & \alpha_4 & - q \alpha_3 \end{array}\right], E = \left[\begin{array}{c} \beta_1 \\ \beta_2 \\ \beta_3 \\ \beta_4 \end{array}\right], \mbox{and} \quad F = \left[\begin{array}{c} 0 \\ 0 \\ 0 \\ 0 \end{array}\right].
\]
We find all $4 \times 4$ minors of D (using Wolfram Mathematica). They are:\\

\[
{\alpha_3}^2 (\alpha_2 \alpha_3 - \alpha_1 \alpha_4) (-1 + q) (1 + q)
\]
\[
\alpha_2 \alpha_3 (\alpha_2 \alpha_3 - \alpha_1 \alpha_4) (-1 + q) (1 + q)
\]
\[
 - \alpha_2 \alpha_4 (\alpha_2 \alpha_3 - \alpha_1 \alpha_4) (-1 + q) (1 + q)
\]
\[
- \alpha_3 \alpha_4 (\alpha_2 \alpha_3 - \alpha_1 \alpha_4) (-1 + q) (1 + q)
\]
\[
\alpha_1 \alpha_2 (\alpha_2 \alpha_3 - \alpha_1 \alpha_4) (-1 + q) q (1 + q)
\]
\[
- {\alpha_2}^2 (\alpha_2 \alpha_3 - \alpha_1 \alpha_4) (-1 + q) (1 + q)
\]
\[
\alpha_1 \alpha_3 (\alpha_2 \alpha_3 - \alpha_1 \alpha_4) (-1 + q) q (1 + q)
\]
Consider the first equation, therefore we have $\alpha_3 = 0$ or $\alpha_2 \alpha_3 - \alpha_1 \alpha_4 = 0$. If $\alpha_3 = 0$, then by Mathematica, we have\\
\[
\alpha_1 {\alpha_2}^2 \alpha_4 (-1 + q) (1 + q), \quad - {\alpha_1}^2 \alpha_2 \alpha_4 (-1 + q) q (1 + q), \quad \alpha_1 \alpha_2 {\alpha_4}^2 (-1 + q) (1 + q).
\]
Therefore in this case, solutions are \\
\[
\{ (0,\beta,0,\delta) \in \mathbb{P}^3 : (\beta,\delta) \in  \mathbb{P}^1 \} \cup \{ (\alpha,0,0,\delta) \in \mathbb{P}^3 : (\alpha,\delta) \in \mathbb{P}^1 \} \cup \{ (\alpha,\beta,0,0) \in \mathbb{P}^3 : (\alpha,\beta) \in \mathbb{P}^1 \}.
\]

\noindent If $\alpha_3 \neq 0$ and $\alpha_2 \alpha_3 = \alpha_1 \alpha_4$, then we obtain \\
\[
\{ (\alpha,\beta,\gamma,\delta) \in \mathbb{P}^3 : \beta \gamma = \alpha \delta \}
\]

\noindent So, in general, the solutions are:\\
\[
\{ (\alpha,0,0,\delta) \in \mathbb{P}^3 : (\alpha,\delta) \in  \mathbb{P}^1 \} \cup \{ (\alpha,\beta,\gamma,\delta) \in \mathbb{P}^3 : \beta \gamma = \alpha \delta \}
\]

\noindent , and the point scheme is $\mathcal V(x_2,x_3) \cup \mathcal V(x_2 x_3 - x_1 x_4)$.
\hfill $\blacksquare$

\subsection{Proposition 2:}

\noindent Suppose $q \in \mathbb{K}^{\times}$. If the algebra
\[
A = \frac{\mathbb{K} \langle x_1,x_2,x_3,x_4 \rangle}{\langle g_1, \ldots, g_6 \rangle}
\]
where
\[
\begin{array}{ll}
g_1 = x_1 x_2 - x_2 x_1, & \quad g_2 = x_3 x_2 - x_2 x_3,\\
g_3 = x_1 x_3 - x_3 x_1, & \quad g_4 = x_4 x_1 - x_1 x_4 + q (x_4 x_3 - x_1 x_2),\\
g_5 = x_4 x_2 - x_2 x_4, & \quad g_6 = x_4 x_3 - x_3 x_4,
\end{array}
\]
, then $A$ has point scheme given by $\mathcal V(x_2 (x_1 x_2 - x_3 x_4), x_3 (x_1 x_2 - x_3 x_4))$ which contains the double line $\mathcal V(x_2,x_3)$ (see Figure 2).

\begin{figure}[ht]
\centering
\includegraphics[width=2cm,height=2cm]{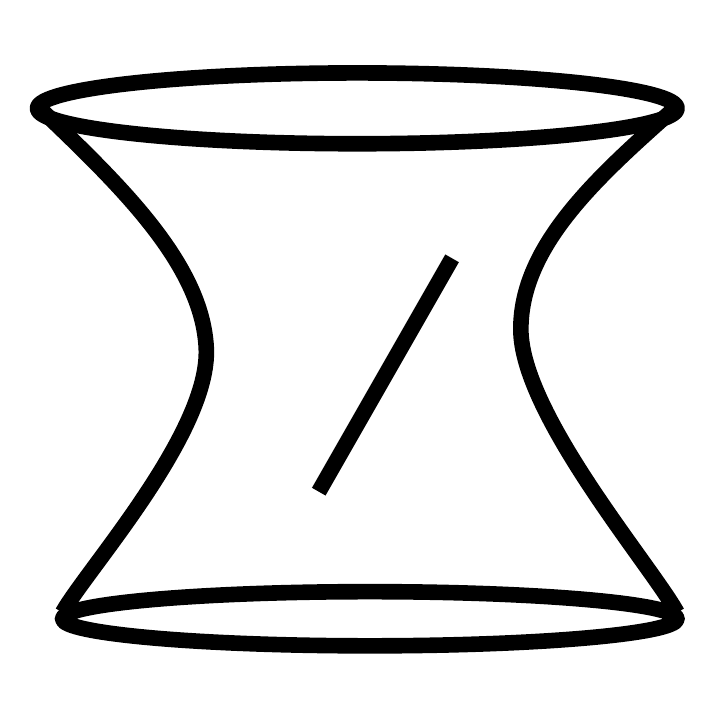}\\
\caption{Depiction of the Point Scheme in Proposition 2}
\end{figure}

\noindent Proof:\\
Suppose
\[
p = ((\alpha_1,\alpha_2,\alpha_3,\alpha_4),(\beta_1,\beta_2,\beta_3,\beta_4)) \in \mathbb{P}^3 \times \mathbb{P}^3.
\]
To find the point scheme $\mathcal{P}$ of $A$, we solve
\[
0 = g_1(p) = \alpha_1 \beta_2 - \alpha_2 \beta_1,
\]
\[
0 = g_2(p) = \alpha_3 \beta_2 - \alpha_2 \beta_3,
\]
\[
0 = g_3(p) = \alpha_1 \beta_3 - \alpha_3 \beta_1,
\]
\[
0 = g_4(p) = \alpha_4 \beta_1 - \alpha_1 \beta_4 + q (\alpha_4 \beta_3 - \alpha_1 \beta_2),
\]
\[
0 = g_5(p) = \alpha_4 \beta_2 - \alpha_2 \beta_4,
\]
\[
0 = g_6(p) = \alpha_4 \beta_3 - \alpha_3 \beta_4,
\]
which yields $D E = F$, where
\[D =
\left[\begin{array}{cccc} - \alpha_2 & \alpha_1 & 0 & 0 \\ 0 &  \alpha_3 & - \alpha_2 & 0 \\ - \alpha_3 & 0 & \alpha_1 & 0 \\ \alpha_4 & - q \alpha_1 & q \alpha_4 & - \alpha_1 \\ 0 & \alpha_4 & 0 & - \alpha_2 \\ 0 & 0 & \alpha_4 & - \alpha_3 \end{array}\right], E = \left[\begin{array}{c} \beta_1 \\ \beta_2 \\ \beta_3 \\ \beta_4 \end{array}\right], \mbox{and} \quad F = \left[\begin{array}{c} 0 \\ 0 \\ 0 \\ 0 \end{array}\right].
\]

\noindent we find all $4 \times 4$ minors of D (using Wolfram Mathematica). They are:\\
  
\[
{\alpha_3}^2 (\alpha_1 \alpha_2 - \alpha_3 \alpha_4) q
\]
\[
- {\alpha_2}^2 (\alpha_1 \alpha_2 - \alpha_3 \alpha_4) q
\]
\[
- \alpha_2 \alpha_3 (\alpha_1 \alpha_2 - \alpha_3 \alpha_4) q
\]
\[\alpha_1 \alpha_2 (\alpha_1 \alpha_2 - \alpha_3 \alpha_4) q
\]
\[
\alpha_1 \alpha_3 (\alpha_1 \alpha_2 - \alpha_3 \alpha_4) q
\]
\[
\alpha_2 \alpha_4 (\alpha_1 \alpha_2 - \alpha_3 \alpha_4) q
\]
\[
\alpha_3 \alpha_4 (\alpha_1 \alpha_2 - \alpha_3 \alpha_4) q
\]
Consider the first equation, therefore we have $\alpha_3 = 0$ or $\alpha_1 \alpha_2 - \alpha_3 \alpha_4 = 0$. If $\alpha_3 = 0$, then by Mathematica, we have\\
\[
- \alpha_1 {\alpha_2}^3 q, \quad {\alpha_1}^2 {\alpha_2}^2 q, \quad \alpha_1 {\alpha_2}^2 \alpha_4 q.
\]
Therefore in this case, solutions are \\
\[
\{ (0,\beta,0,\delta) \in \mathbb{P}^3 : (\beta,\delta) \in  \mathbb{P}^1 \} \cup \{ (\alpha,0,0,\delta) \in \mathbb{P}^3 : (\alpha,\delta) \in \mathbb{P}^1 \}.
\]

\noindent If $\alpha_3 \neq 0$ and $\alpha_1 \alpha_2 = \alpha_3 \alpha_4$, then we obtain \\
\[
\{ (\alpha,\beta,\gamma,\delta) \in \mathbb{P}^3 : \alpha \beta = \gamma \delta \}
\]

\noindent So, in general, the solutions are:\\
\[
\{ (\alpha,\beta,\gamma,\delta) \in \mathbb{P}^3 : \alpha \beta = \gamma \delta \}
\]

\noindent , and the point scheme is $\mathcal V(x_2 (x_1 x_2 - x_3 x_4), x_3 (x_1 x_2 - x_3 x_4))$ which contains the double line $\mathcal V(x_2,x_3)$.\\

\noindent Notice that
\[
\langle x_1 x_2 - x_3 x_4 \rangle \langle {x_2}^2,x_2 x_3,x_1 x_2,x_1 x_3, x_2 x_4, {x_3}^2,x_3 x_4 \rangle
\]
\[
= \langle x_1 x_2 - x_3 x_4 \rangle \langle x_2, x_3 \rangle \langle x_1,x_2,x_3,x_4 \rangle
\]
Therefore,  the point scheme $\mathcal{P}$ of $A$ is\\
\[
\mathcal V(x_1 x_2 - x_3 x_4) \cup \mathcal V(x_2, x_3) \cup \mathcal V(x_1, x_2, x_3, x_4)
\]
\[
= \mathcal V(x_1 x_2 - x_3 x_4) \cup \mathcal V(x_2, x_3)
\]
\[
\mbox{where} \qquad \mathcal V(x_2, x_3) \subset \mathcal V(x_1 x_2 - x_3 x_4).
\]

\hfill $\blacksquare$

\subsection{Proposition 3:}

\noindent Suppose $q \in \mathbb{K}^{\times} \setminus \{ -1 \}$. If the algebra
\[
A = \frac{\mathbb{K} \langle x_1,x_2,x_3,x_4 \rangle}{\langle g_1, \ldots, g_6 \rangle}
\]
where
\[
\begin{array}{ll}
g_1 = x_1 x_2 - x_2 x_1, & \quad g_2 = x_2 x_3 - x_3 x_2,\\
g_3 = x_1 x_3 - x_3 x_1, & \quad g_4 = x_1 x_4 - x_4 x_1,\\
g_5 = x_2 x_4 - x_4 x_2 - q ({x_1}^2 - x_2 x_4), & \quad g_6 = x_4 x_3 - x_3 x_4,
\end{array}
\]
, then $A$ has point scheme given by $Q \cup L$ where $Q = \mathcal V({x_1}^2 - x_2 x_4)$ and $L = \mathcal V(x_1,x_3)$ (see Figure 3).

\begin{figure}[ht]
\centering
\includegraphics[width=2cm,height=2cm]{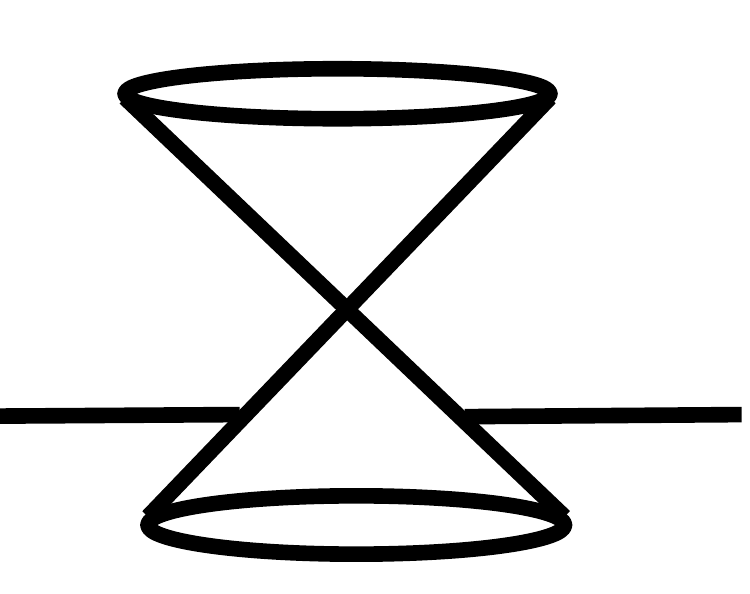}\\
\caption{Depiction of the Point Scheme in Proposition 3}
\end{figure}

\noindent Proof:\\

\noindent Suppose
\[
p = ((\alpha_1,\alpha_2,\alpha_3,\alpha_4),(\beta_1,\beta_2,\beta_3,\beta_4)) \in \mathbb{P}^3 \times \mathbb{P}^3.
\]
To find the point scheme $\mathcal{P}$ of $A$, we solve
\[
0 = g_1(p) = \alpha_1 \beta_2 - \alpha_2 \beta_1,
\]
\[
0 = g_2(p) = \alpha_2 \beta_3 - \alpha_3 \beta_2,
\]
\[
0 = g_3(p) = \alpha_1 \beta_3 - \alpha_3 \beta_1,
\]
\[
0 = g_4(p) = \alpha_1 \beta_4 - \alpha_4 \beta_1,
\]
\[
0 = g_5(p) = \alpha_2 \beta_4 - \alpha_4 \beta_2 - q (\alpha_1 \beta_1 - \alpha_2 \beta_4),
\]
\[
0 = g_6(p) = \alpha_4 \beta_3 - \alpha_3 \beta_4,
\]
which yields $D E = F$, where
\[D =
\left[\begin{array}{cccc} - \alpha_2 & \alpha_1 & 0 & 0 \\ 0 & - \alpha_3 & \alpha_2 & 0 \\ - \alpha_3 & 0 & \alpha_1 & 0 \\ - \alpha_4 & 0 & 0 & \alpha_1 \\ - q \alpha_1 & - \alpha_4 & 0 & (q + 1) \alpha_2 \\ 0 & 0 & \alpha_4 & - \alpha_3 \end{array}\right], E = \left[\begin{array}{c} \beta_1 \\ \beta_2 \\ \beta_3 \\ \beta_4 \end{array}\right], \mbox{and} \quad F = \left[\begin{array}{c} 0 \\ 0 \\ 0 \\ 0 \end{array}\right].
\]

\noindent We find all $4 \times 4$ minors of D (using Wolfram Mathematica). They are:\\
 
\[
{\alpha_1}^2 ({\alpha_1}^2 - \alpha_2 \alpha_4) q
\]
\[
- {\alpha_3}^2 ({\alpha_1}^2 - \alpha_2 \alpha_4) q
\]
\[
\alpha_1 \alpha_2 ({\alpha_1}^2 - \alpha_2 \alpha_4) q
\]
\[
\alpha_2 \alpha_3 ({\alpha_1}^2 - \alpha_2 \alpha_4) q
\]
\[
\alpha_1 \alpha_3 ({\alpha_1}^2 - \alpha_2 \alpha_4) q
\]
\[
\alpha_1 \alpha_4 ({\alpha_1}^2 - \alpha_2 \alpha_4) q
\]
\[
- \alpha_3 \alpha_4 ({\alpha_1}^2 - \alpha_2 \alpha_4) q
\]

\vspace{0.2in}

\noindent Consider the first equation, therefore we have $\alpha_1 = 0$ or ${\alpha_1}^2 - \alpha_2 \alpha_4 = 0$. If $\alpha_1 = 0$, then by Mathematica, we have\\
\[
- {\alpha_2}^2 \alpha_3 \alpha_4 q, \quad \alpha_2 {\alpha_3}^2 \alpha_4 q, \quad \alpha_2 \alpha_3 {\alpha_4}^2 q.
\]
\noindent Therefore in this case, solutions are \\
\[
\{ (0,0,\gamma,\delta) \in \mathbb{P}^3 : (\gamma,\delta) \in  \mathbb{P}^1 \} \cup \{ (0,\beta,0,\delta) \in \mathbb{P}^3 : (\beta,\delta) \in \mathbb{P}^1 \} \cup \{ (0,\beta,\gamma,0) \in \mathbb{P}^3 : (\beta,\gamma) \in \mathbb{P}^1 \}.
\]

\noindent If $\alpha_3 \neq 0$ and ${\alpha_1}^2 = \alpha_2 \alpha_4$, then we obtain \\
\[
\{ (\alpha,\beta,\gamma,\delta) \in \mathbb{P}^3 : \alpha^2 = \beta \delta \}
\]

\noindent So, in general, the solutions are:\\
\[
\{ (0,\beta,0,\delta) \in \mathbb{P}^3 : (\beta,\delta) \in  \mathbb{P}^1 \} \cup \{ (\alpha,\beta,\gamma,\delta) \in \mathbb{P}^3 : \alpha^2 = \beta \delta \}
\]

\noindent , and the point scheme is $\mathcal V({x_1}^2 - x_2 x_4) \cup \mathcal V(x_1,x_3)$
\hfill $\blacksquare$

\subsection{Proposition 4:}

\noindent If the algebra
\[
A = \frac{\mathbb{K} \langle x_1,x_2,x_3,x_4 \rangle}{\langle g_1, \ldots, g_6 \rangle}
\]
where
\[
\begin{array}{ll}
g_1 = x_1 x_2 - x_2 x_1, & \quad g_2 = x_3 x_2 - x_2 x_3,\\
g_3 = x_1 x_3 - x_3 x_1, & \quad g_4 = x_1 x_4 - x_4 x_1 - {x_1}^2 + x_4 x_3,\\
g_5 = x_2 x_4 - x_4 x_2, & \quad g_6 = x_3 x_4 - x_4 x_3,
\end{array}
\]
, then $A$ has point scheme given by $Q \cup L$ where $Q = \mathcal V({x_1}^2 - x_3 x_4)$ and $L = \mathcal V(x_2,x_3)$ (so the line $L$ is tangential to the quadric $Q$ at a nonsingular point of $Q$)(see Figure 4).

\begin{figure}[ht]
\centering
\includegraphics[width=2cm,height=2cm]{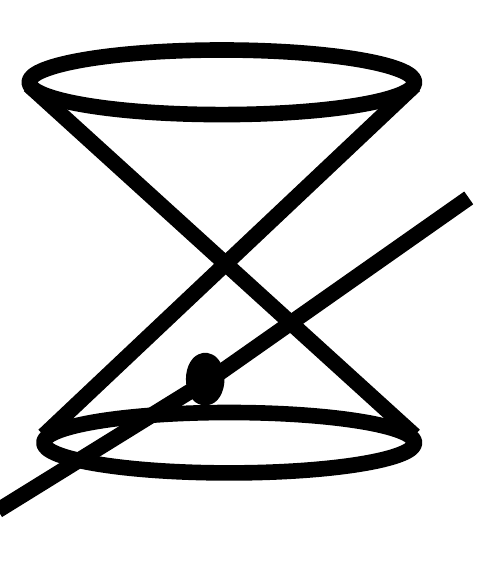}\\
\caption{Depiction of the Point Scheme in Proposition 4}
\end{figure}

\noindent Proof:\\

\noindent Suppose
\[
p = ((\alpha_1,\alpha_2,\alpha_3,\alpha_4),(\beta_1,\beta_2,\beta_3,\beta_4)) \in \mathbb{P}^3 \times \mathbb{P}^3.
\]
To find the point scheme $\mathcal{P}$ of $A$, we solve
\[
0 = g_1(p) = \alpha_1 \beta_2 - \alpha_2 \beta_1,
\]
\[
0 = g_2(p) = \alpha_3 \beta_2 - \alpha_2 \beta_3,
\]
\[
0 = g_3(p) = \alpha_1 \beta_3 - \alpha_3 \beta_1,
\]
\[
0 = g_4(p) = \alpha_1 \beta_4 - \alpha_4 \beta_1 - \alpha_1 \beta_1 + \alpha_4 \beta_3,
\]
\[
0 = g_5(p) = \alpha_2 \beta_4 - \alpha_4 \beta_2,
\]
\[
0 = g_6(p) = \alpha_3 \beta_4 - \alpha_4 \beta_3,
\]
which yields $D E = F$, where
\[D =
\left[\begin{array}{cccc} - \alpha_2 & \alpha_1 & 0 & 0 \\ 0 & \alpha_3 & - \alpha_2 & 0 \\ - \alpha_3 & 0 & \alpha_1 & 0 \\ - \alpha_4 - \alpha_1 & 0 & \alpha_4 & \alpha_1 \\ 0 & - \alpha_4 & 0 & \alpha_2 \\ 0 & 0 & - \alpha_4 & \alpha_3 \end{array}\right], E = \left[\begin{array}{c} \beta_1 \\ \beta_2 \\ \beta_3 \\ \beta_4 \end{array}\right], \mbox{and} \quad F = \left[\begin{array}{c} 0 \\ 0 \\ 0 \\ 0 \end{array}\right].
\]

\noindent We find all $4 \times 4$ minors of D (using Wolfram Mathematica). They are:\\
 
\[
- {\alpha_3}^2 ({\alpha_1}^2 - \alpha_3 \alpha_4)
\]
\[
{\alpha_2}^2 ({\alpha_1}^2 - \alpha_3 \alpha_4)
\]
\[
\alpha_2 \alpha_3 ({\alpha_1}^2 - \alpha_3 \alpha_4)
\]
\[
 - \alpha_1 \alpha_2 ({\alpha_1}^2 - \alpha_3 \alpha_4)
\]
\[
- \alpha_1 \alpha_3 ({\alpha_1}^2 - \alpha_3 \alpha_4)
\]
\[
\alpha_2 \alpha_4 ({\alpha_1}^2 - \alpha_3 \alpha_4)
\]
\[
\alpha_3 \alpha_4 ({\alpha_1}^2 - \alpha_3 \alpha_4)
\]

\vspace{0.2in}

\noindent Consider the first equation, therefore we have $\alpha_3 = 0$ or ${\alpha_1}^2 - \alpha_3 \alpha_4 = 0$. If $\alpha_3 = 0$, then by Mathematica, we have\\
\[
{\alpha_1}^2 {\alpha_2}^2, \quad - {\alpha_1}^3 \alpha_2, \quad {\alpha_1}^2 \alpha_2 \alpha_4.
\]
\noindent Therefore in this case, solutions are \\
\[
\{ (0,\beta,0,\delta) \in \mathbb{P}^3 : (\beta,\delta) \in  \mathbb{P}^1 \} \cup \{ (\alpha,0,0,\delta) \in \mathbb{P}^3 : (\alpha,\delta) \in \mathbb{P}^1 \}.
\]

\noindent If $\alpha_3 \neq 0$ and ${\alpha_1}^2 = \alpha_3 \alpha_4$, then we obtain \\
\[
\{ (\alpha,\beta,\gamma,\delta) \in \mathbb{P}^3 : {\alpha}^2 = \gamma \delta \}
\]

\noindent So, in general, the solutions are:\\
\[
\{ (\alpha,0,0,\delta) \in \mathbb{P}^3 : (\alpha,\delta) \in  \mathbb{P}^1 \} \cup \{ (\alpha,\beta,\gamma,\delta) \in \mathbb{P}^3 : \alpha^2 = \gamma \delta \}
\]

\noindent , and the point scheme is $\mathcal V({x_1}^2 - x_3 x_4) \cup \mathcal V(x_2,x_3)$.
\hfill $\blacksquare$

\subsection{Proposition 5:}

\noindent If the algebra
\[
A = \frac{\mathbb{K} \langle x_1,x_2,x_3,x_4 \rangle}{\langle g_1, \ldots, g_6 \rangle}
\]
where
\[
\begin{array}{ll}
g_1 = x_1 x_2 - x_2 x_1, & \quad g_2 = x_2 x_3 - x_3 x_2,\\
g_3 = x_1 x_3 - x_3 x_1, & \quad g_4 = x_1 x_4 - x_4 x_1 - {x_1}^2 + x_2 x_3,\\
g_5 = x_2 x_4 - x_4 x_2, & \quad g_6 = x_3 x_4 - x_4 x_3,
\end{array}
\]
, then $A$ has point scheme given by $Q \cup L$ where $Q = \mathcal V({x_1}^2 - x_2 x_3)$ and $L = \mathcal V(x_2,x_3)$ (so the line $L$ is tangential to the quadric $Q$ at a singular point of $Q$)(see Figure 5).

\begin{figure}[ht]
\centering
\includegraphics[width=2cm,height=2cm]{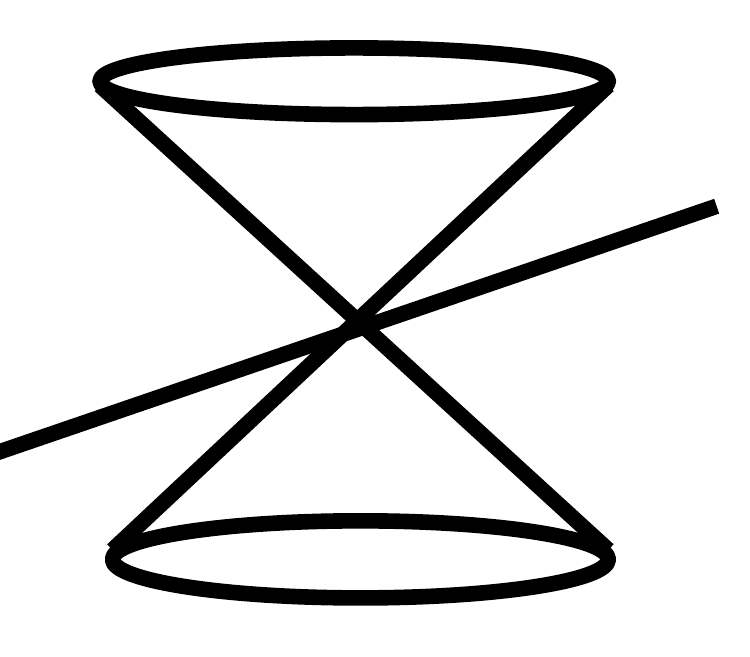}\\
\caption{Depiction of the Point Scheme in Proposition 5}
\end{figure}

\noindent Proof:\\

\noindent Suppose
\[
p = ((\alpha_1,\alpha_2,\alpha_3,\alpha_4),(\beta_1,\beta_2,\beta_3,\beta_4)) \in \mathbb{P}^3 \times \mathbb{P}^3.
\]
To find the point scheme $\mathcal{P}$ of $A$, we solve
\[
0 = g_1(p) = \alpha_1 \beta_2 - \alpha_2 \beta_1,
\]
\[
0 = g_2(p) = \alpha_2 \beta_3 - \alpha_3 \beta_2,
\]
\[
0 = g_3(p) = \alpha_1 \beta_3 - \alpha_3 \beta_1,
\]
\[
0 = g_4(p) = \alpha_1 \beta_4 - \alpha_4 \beta_1 - \alpha_1 \beta_1 + \alpha_2 \beta_3,
\]
\[
0 = g_5(p) = \alpha_2 \beta_4 - \alpha_4 \beta_2,
\]
\[
0 = g_6(p) = \alpha_3 \beta_4 - \alpha_4 \beta_3,
\]
which yields $D E = F$, where
\[D =
\left[\begin{array}{cccc} - \alpha_2 & \alpha_1 & 0 & 0 \\ 0 & - \alpha_3 & \alpha_2 & 0 \\ - \alpha_3 & 0 & \alpha_1 & 0 \\ - \alpha_4 - \alpha_1 & 0 & \alpha_2 & \alpha_1 \\ 0 & - \alpha_4 & 0 & \alpha_2 \\ 0 & 0 & - \alpha_4 & \alpha_3 \end{array}\right], E = \left[\begin{array}{c} \beta_1 \\ \beta_2 \\ \beta_3 \\ \beta_4 \end{array}\right], \mbox{and} \quad F = \left[\begin{array}{c} 0 \\ 0 \\ 0 \\ 0 \end{array}\right].
\]

\noindent We find all $4 \times 4$ minors of D (using Wolfram Mathematica). They are:\\

\[
{\alpha_3}^2 ({\alpha_1}^2 - \alpha_2 \alpha_3)
\]
\[
- {\alpha_2}^2 ({\alpha_1}^2 - \alpha_2 \alpha_3)
\]
\[
\alpha_2 \alpha_3 ({\alpha_1}^2 - \alpha_2 \alpha_3)
\]
\[
- \alpha_1 \alpha_2 ({\alpha_1}^2 - \alpha_2 \alpha_3)
\]
\[
- \alpha_1 \alpha_3 ({\alpha_1}^2 - \alpha_2 \alpha_3)
\]
\[
\alpha_2 \alpha_4 ({\alpha_1}^2 - \alpha_2 \alpha_3)
\]
\[
\alpha_3 \alpha_4 ({\alpha_1}^2 - \alpha_2 \alpha_3)
\]

\vspace{0.2in}

\noindent Consider the first equation, therefore we have $\alpha_3 = 0$ or ${\alpha_1}^2 - \alpha_2 \alpha_3 = 0$. If $\alpha_3 = 0$, then by Mathematica, we have\\
\[
- {\alpha_1}^2 {\alpha_2}^2, \quad - {\alpha_1}^3 \alpha_2, \quad {\alpha_1}^2 \alpha_2 \alpha_4.
\]
\noindent Therefore in this case, solutions are \\
\[
\{ (0,\beta,0,\delta) \in \mathbb{P}^3 : (\beta,\delta) \in  \mathbb{P}^1 \} \cup \{ (\alpha,0,0,\delta) \in \mathbb{P}^3 : (\alpha,\delta) \in \mathbb{P}^1 \}.
\]

\noindent If $\alpha_3 \neq 0$ and ${\alpha_1}^2 = \alpha_2 \alpha_3$, then we obtain \\
\[
\{ (\alpha,\beta,\gamma,\delta) \in \mathbb{P}^3 : \alpha^2 = \beta \gamma \}
\]

\noindent So, in general, the solutions are:\\
\[
\{ (\alpha,0,0,\delta) \in \mathbb{P}^3 : (\alpha,\delta) \in  \mathbb{P}^1 \} \cup \{ (\alpha,\beta,\gamma,\delta) \in \mathbb{P}^3 : \alpha^2 = \beta \gamma \}
\]

\noindent , and the point scheme is $\mathcal V({x_1}^2 - x_2 x_3) \cup \mathcal V(x_2,x_3)$.
\hfill $\blacksquare$

\subsection{Proposition 6:}

\noindent If the algebra
\[
A = \frac{\mathbb{K} \langle x_1,x_2,x_3,x_4 \rangle}{\langle g_1, \ldots, g_6 \rangle}
\]
where
\[
\begin{array}{ll}
g_1 = x_1 x_2 - x_2 x_1, & \quad g_2 = x_2 x_3 - x_3 x_2,\\
g_3 =  x_1 x_3 - x_3 x_1, & \quad g_4 = x_1 x_4 - x_4 x_1,\\
g_5 = x_2 x_4 - x_4 x_2, & \quad g_6 = x_3 x_4 - x_4 x_3 - {x_1}^2 + x_2 x_3,
\end{array}
\]
, then $A$ has point scheme given by $\mathcal V (x_1 ({x_1}^2 - x_2 x_3), x_2 ({x_1}^2 - x_2 x_3))$, which contains the double line $V(x_1,x_2)$ (see Figure 6).

\begin{figure}[ht]
\centering
\includegraphics[width=2cm,height=2cm]{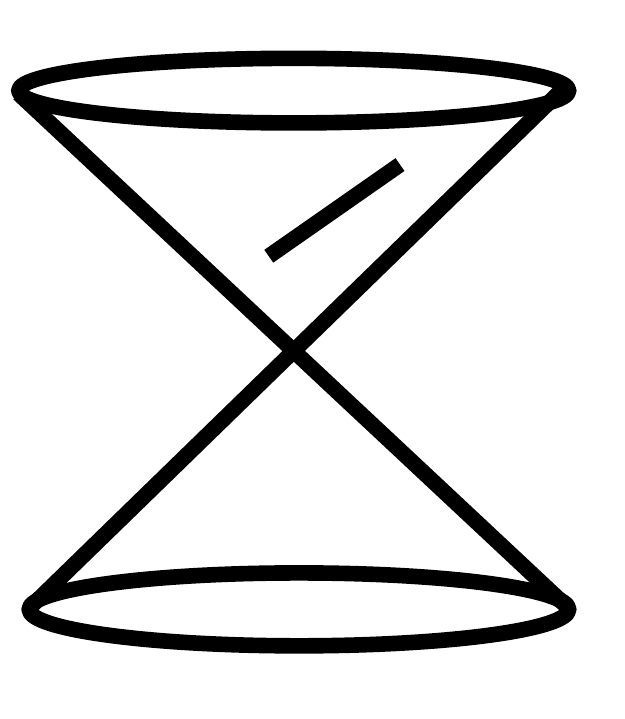}\\
\caption{Depiction of the Point Scheme in Proposition 6}
\end{figure}

\noindent Proof:\\

\noindent Suppose
\[
p = ((\alpha_1,\alpha_2,\alpha_3,\alpha_4),(\beta_1,\beta_2,\beta_3,\beta_4)) \in \mathbb{P}^3 \times \mathbb{P}^3.
\]
To find the point scheme $\mathcal{P}$ of $A$, we solve
\[
0 = g_1(p) = \alpha_1 \beta_2 - \alpha_2 \beta_1,
\]
\[
0 = g_2(p) = \alpha_2 \beta_3 - \alpha_3 \beta_2,
\]
\[
0 = g_3(p) = \alpha_1 \beta_3 - \alpha_3 \beta_1,
\]
\[
0 = g_4(p) = \alpha_1 \beta_4 - \alpha_4 \beta_1,
\]
\[
0 = g_5(p) = \alpha_2 \beta_4 - \alpha_4 \beta_2,
\]
\[
0 = g_6(p) = \alpha_3 \beta_4 - \alpha_4 \beta_3 - \alpha_1 \beta_1 + \alpha_2 \beta_3,
\]
which yields $D E = F$, where
\[D =
\left[\begin{array}{cccc} - \alpha_2 & \alpha_1 & 0 & 0 \\ 0 & - \alpha_3 & \alpha_2 & 0 \\ - \alpha_3 & 0 & \alpha_1 & 0 \\ - \alpha_4 & 0 & 0 & \alpha_1 \\ 0 & - \alpha_4 & 0 & \alpha_2 \\ - \alpha_1 & 0 & \alpha_2 - \alpha_4 & \alpha_3 \end{array}\right], E = \left[\begin{array}{c} \beta_1 \\ \beta_2 \\ \beta_3 \\ \beta_4 \end{array}\right], \mbox{and} \quad F = \left[\begin{array}{c} 0 \\ 0 \\ 0 \\ 0 \end{array}\right].
\]

\noindent We find all $4 \times 4$ minors of D (using Wolfram Mathematica). They are:\\

\[
{\alpha_1}^2 ({\alpha_1}^2 - \alpha_2 \alpha_3)
\]
\[
\alpha_1 \alpha_2 ({\alpha_1}^2 - \alpha_2 \alpha_3)
\]
\[
{\alpha_2}^2 ({\alpha_1}^2 - \alpha_2 \alpha_3)
\]
\[
- \alpha_1 \alpha_3 ({\alpha_1}^2 - \alpha_2 \alpha_3)
\]
\[
- \alpha_2 \alpha_3 ({\alpha_1}^2 - \alpha_2 \alpha_3)
\]
\[
- \alpha_2 \alpha_4({\alpha_1}^2 - \alpha_2 \alpha_3)
\]
\[
- \alpha_1 \alpha_4 ({\alpha_1}^2 - \alpha_2 \alpha_3)
\]

\vspace{0.2in}

\noindent Consider the first equation, therefore we have $\alpha_1 = 0$ or ${\alpha_1}^2 - \alpha_2 \alpha_3 = 0$. If $\alpha_1 = 0$, then by Mathematica, we have\\
\[
- {\alpha_2}^3 \alpha_3, \quad {\alpha_2}^2 {\alpha_3}^2, \quad {\alpha_2}^2 \alpha_3 \alpha_4.
\]
\noindent Therefore in this case, solutions are \\
\[
\{ (0,0,\gamma,\delta) \in \mathbb{P}^3 : (\gamma,\delta) \in \mathbb{P}^1 \} \cup \{ (0,\beta,0,\delta) \in \mathbb{P}^3 : (\beta,\delta) \in  \mathbb{P}^1 \}.
\]

\noindent If $\alpha_1 \neq 0$ and ${\alpha_1}^2 = \alpha_2 \alpha_3$, then we obtain \\
\[
\{ (\alpha,\beta,\gamma,\delta) \in \mathbb{P}^3 : \alpha^2 = \beta \gamma \}
\]

\noindent So, in general, the solutions are:\\
\[
\{ (\alpha,\beta,\gamma,\delta) \in \mathbb{P}^3 : \alpha^2 = \beta \gamma \}
\]

\noindent , and the point scheme is $\mathcal V (x_1 ({x_1}^2 - x_2 x_3), x_2 ({x_1}^2 - x_2 x_3))$, which contains the double line $\mathcal V(x_1,x_2)$.\\

\noindent Notice that
\[
\langle {x_1}^2 - x_2 x_3 \rangle \langle x_1 x_2,{x_2}^2,{x_1}^2,x_1 x_3,x_2 x_3, x_2 x_4,x_1 x_4 \rangle
\]
\[
= \langle {x_1}^2 - x_2 x_3 \rangle \langle x_1, x_2 \rangle \langle x_1,x_2,x_3,x_4 \rangle
\]
Therefore,  the point scheme $\mathcal{P}$ of $A$ is\\
\[
\mathcal V({x_1}^2 - x_2 x_3) \cup \mathcal V(x_1, x_2) \cup \mathcal V(x_1, x_2, x_3, x_4)
\]
\[
= \mathcal V({x_1}^2 - x_2 x_3) \cup \mathcal V(x_1, x_2)
\]
\[
\mbox{where} \qquad \mathcal V(x_1, x_2) \subset \mathcal V({x_1}^2 - x_2 x_3).
\]

\hfill $\blacksquare$

\pagebreak

\section{Appendix}

Mathematica code for Proposition 1.

\begin{figure}[ht]
\includegraphics[width=14cm,height=14cm]{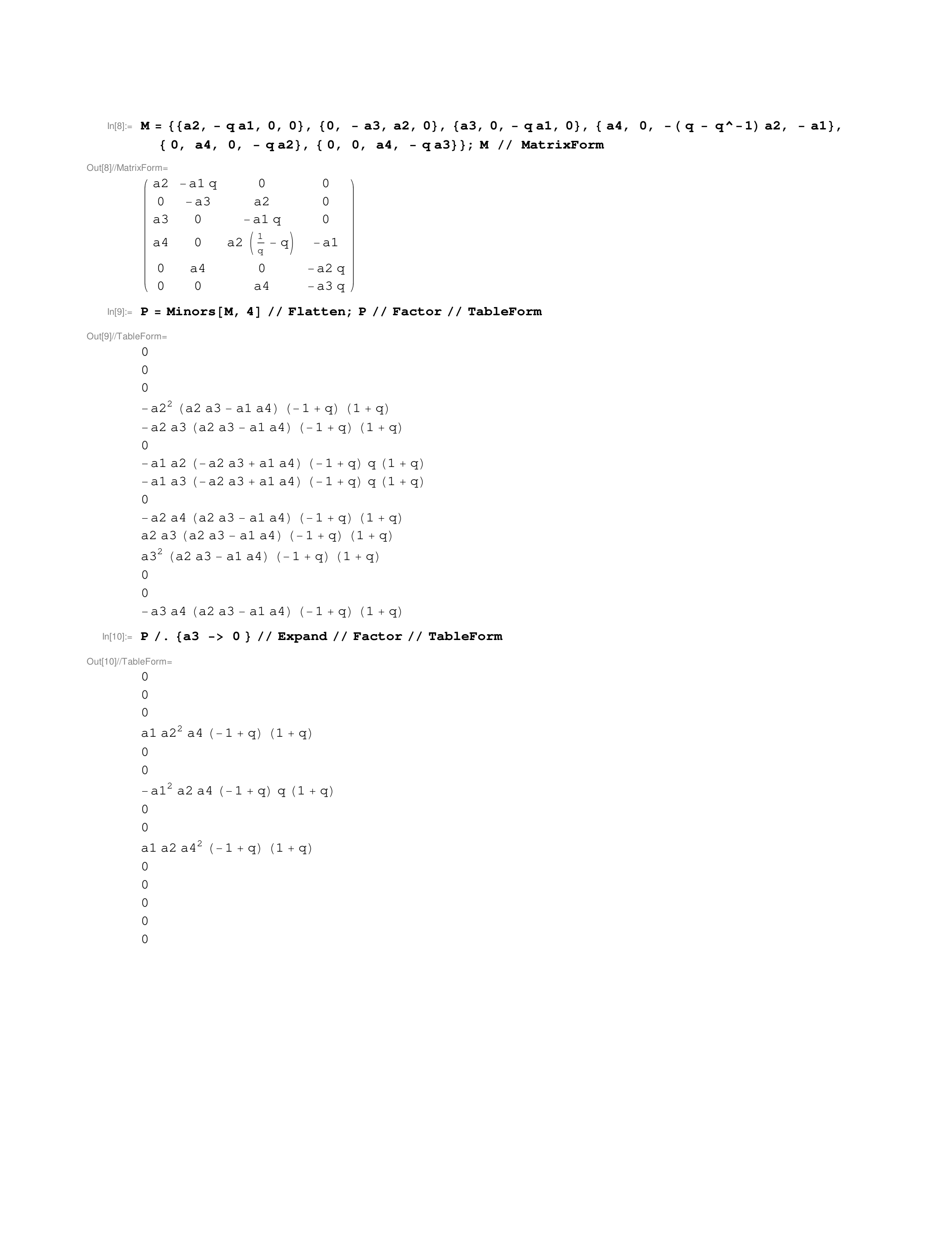}\\
\end{figure}

% Generate the bibliography.
%\bibliography{latex-sample}
%\bibliographystyle{unsrt}

% Generate the bibliography.
%\bibliography{latex-sample}
%\bibliographystyle{unsrt}

\pagebreak

\end{document}